\documentclass[a4paper,12pt]{article}
\usepackage[all]{xy}
\usepackage[T1]{fontenc}
\usepackage[dvips]{graphicx}
\usepackage{amsmath,amssymb,amsfonts,amsthm,latexsym,mathrsfs,textcomp,verbatim}
\xyoption{web}

\title{A characterization theorem for geometric logic}
\author{{Olivia Caramello} \vspace{3 mm}\\ {\small Centro di Ricerca Matematica Ennio De Giorgi}\\{\small Scuola Normale Superiore} \\{\small Piazza dei Cavalieri, 3, 56100 Pisa, Italy}\\{\small olivia.caramello@sns.it}}
\date{\today}
\begin{document}
\bgroup           % fake a titlepage 
\let\footnoterule\relax  % no rule above thanks footnotes 
\maketitle

\begin{abstract}
We establish a criterion for deciding whether a class of structures is the class of models of a geometric theory inside Grothendieck toposes; then we specialize this result to obtain a characterization of the infinitary first-order theories which are geometric in terms of their models in Grothendieck toposes, solving a problem posed by Ieke Moerdijk in 1989.      
\end{abstract} 
\egroup 

\vspace{5 mm}

%MACROS-----------------------------------------------------------------------------------------------------------------------

%	European dates ``19 April 1990'' not ``April 19, 1990''
\def\Monthnameof#1{\ifcase#1\or
   January\or February\or March\or April\or May\or June\or
   July\or August\or September\or October\or November\or December\fi}
\def\today{\number\day~\Monthnameof\month~\number\year}

%===========================================================================
%	END OF PROOF BOX
%
%
%  The complexity of the macro necessary to get a little box on the
%  right-hand-side at the end of a proof is amazing.  It really does
%  have to be this long!  Otherwise you're liable to get it at the
%  beginning of the next line, or even on the next page.
%
\def\pushright#1{{%        set up
   \parfillskip=0pt            % so \par doesnt push \square to left
   \widowpenalty=10000         % so we dont break the page before \square
   \displaywidowpenalty=10000  % ditto
   \finalhyphendemerits=0      % TeXbook exercise 14.32
  %
  %                 horizontal
   \leavevmode                 % \nobreak means lines not pages
   \unskip                     % remove previous space or glue
   \nobreak                    % don't break lines
   \hfil                       % ragged right if we spill over
   \penalty50                  % discouragement to do so
   \hskip.2em                  % ensure some space
   \null                       % anchor following \hfill
   \hfill                      % push \square to right
   {#1}                        % the end-of-proof mark (or whatever)
  %
  %                   vertical
   \par}}                      % build paragraph

% prefer proofs with statements, also space after
\def\qed{\pushright{$\square$}\penalty-700 \smallskip}

\newtheorem{theorem}{Theorem}[section]

\newtheorem{proposition}[theorem]{Proposition}

\newtheorem{scholium}[theorem]{Scholium}

\newtheorem{lemma}[theorem]{Lemma}

\newtheorem{corollary}[theorem]{Corollary}

\newtheorem{conjecture}[theorem]{Conjecture}

\newenvironment{proofs}%
 {\begin{trivlist}\item[]{\bf Proof }}%
 {\qed\end{trivlist}}

  \newtheorem{rmk}[theorem]{Remark}
\newenvironment{remark}{\begin{rmk}\em}{\end{rmk}}

  \newtheorem{rmks}[theorem]{Remarks}
\newenvironment{remarks}{\begin{rmks}\em}{\end{rmks}}

  \newtheorem{defn}[theorem]{Definition}
\newenvironment{definition}{\begin{defn}\em}{\end{defn}}

  \newtheorem{eg}[theorem]{Example}
\newenvironment{example}{\begin{eg}\em}{\end{eg}}

  \newtheorem{egs}[theorem]{Examples}
\newenvironment{examples}{\begin{egs}\em}{\end{egs}}

%%%%%%%%%%%%%%%%%%%%%%%%%%%%%%%%%%%%%%%%%%%%%%%%%%%%%%%%%%%%%%%%%%%%%%

%  change some single-character symbols to be more appropriate for logic
\mathcode`\<="4268  % < = \langle
\mathcode`\>="5269  % > = \rangle
\mathcode`\.="313A  % make . a binary  relation
\mathchardef\semicolon="603B % the original
\mathchardef\gt="313E
\mathchardef\lt="313C

\newcommand{\app}% application
 {{\sf app}}

\newcommand{\Ass}% category of assemblies
 {{\bf Ass}}

\newcommand{\ASS}% indexed version of \Ass
 {{\mathbb A}{\sf ss}}

\newcommand{\Bb}%blackboard bold
{\mathbb}

\newcommand{\biimp}% bi-implication
 {\!\Leftrightarrow\!}

\newcommand{\bim}% bimorphism
 {\rightarrowtail\kern-1em\twoheadrightarrow}

\newcommand{\bjg}% bi-judgement
 {\mathrel{{\dashv}\,{\vdash}}}

\newcommand{\bstp}[3]% bimorphism version of \stp
 {\mbox{$#1\! : #2 \bim #3$}}

\newcommand{\cat}%concatenation
 {\!\mbox{\t{\ }}}

\newcommand{\cinf}%C-infinity
 {C^{\infty}}

\newcommand{\cinfrg}%category of C-infinity rings
 {\cinf\hy{\bf Rng}}

\newcommand{\cocomma}[2]% cocomma object
 {\mbox{$(#1\!\uparrow\!#2)$}}

\newcommand{\cod}% codomain
 {{\rm cod}}

\newcommand{\comma}[2]% comma object
 {\mbox{$(#1\!\downarrow\!#2)$}}

\newcommand{\comp}% composition
 {\circ}

\newcommand{\cons}% concatenation
 {{\sf cons}}

\newcommand{\Cont}% category of continuous G-sets
 {{\bf Cont}}

\newcommand{\ContE}% continuous G-sets relative to $\cal E$
 {{\bf Cont}_{\cal E}}

\newcommand{\ContS}% ditto for $\cal S$
 {{\bf Cont}_{\cal S}}

\newcommand{\cover}% cover
 {-\!\!\triangleright\,}

\newcommand{\cstp}[3]% cover version of \stp
 {\mbox{$#1\! : #2 \cover #3$}}

\newcommand{\Dec}% decalage
 {{\rm Dec}}

\newcommand{\DEC}% decalage (\Bbb version)
 {{\mathbb D}{\sf ec}}

\newcommand{\den}[1]% denotation of #1 
 {[\![#1]\!]}

\newcommand{\Desc}% category of descent data
 {{\bf Desc}}

\newcommand{\dom}% domain
 {{\rm dom}}

\newcommand{\Eff}% effective topos
 {{\bf Eff}}

\newcommand{\EFF}% indexed version of \Eff
 {{\mathbb E}{\sf ff}}

\newcommand{\empstg}% empty string
 {[\,]}

\newcommand{\epi}% epimorphism
 {\twoheadrightarrow}

\newcommand{\estp}[3]% epimorphism version of \stp
 {\mbox{$#1 \! : #2 \epi #3$}}

\newcommand{\ev}% evaluation
 {{\rm ev}}

\newcommand{\Ext}% category of extracts
 {{\rm Ext}}

\newcommand{\fr}% Fraktur (i.e. Gothic)
 {\sf}

\newcommand{\fst}% first projection
 {{\sf fst}}

\newcommand{\fun}[2]% function-type
 {\mbox{$[#1\!\to\!#2]$}}

\newcommand{\funs}[2]% function-type as subscript
 {[#1\!\to\!#2]}

\newcommand{\Gl}% topos obtained by glueing
 {{\bf Gl}}

\newcommand{\hash}% hash sign (used as infix)
 {\,\#\,}

\newcommand{\hy}% hyphen (in math mode)
 {\mbox{-}}

\newcommand{\im}% image
 {{\rm im}}

\newcommand{\imp}% implication
 {\!\Rightarrow\!}

\newcommand{\Ind}[1]% ind-completion of #1
 {{\rm Ind}\hy #1}

\newcommand{\iten}[1]% enumerated item
{\item[{\rm (#1)}]}

\newcommand{\iter}% iterator
 {{\sf iter}}

\newcommand{\Kalg}%category of $K$-algebras
 {K\hy{\bf Alg}}

\newcommand{\llim}% left (inverse) limit
 {{\mbox{$\lower.95ex\hbox{{\rm lim}}$}\atop{\scriptstyle
{\leftarrow}}}{}}

\newcommand{\llimd}% \displaymath version of \llim
 {\lower0.37ex\hbox{$\pile{\lim \\ {\scriptstyle
\leftarrow}}$}{}}

\newcommand{\Mf}%category of manifolds
 {{\bf Mf}}

\newcommand{\Mod}% category of modest assemblies
 {{\bf Mod}}

\newcommand{\MOD}% indexed version of \Mod
{{\mathbb M}{\sf od}}

\newcommand{\mono}% monomorphism 
 {\rightarrowtail}

\newcommand{\mor}% class of morphisms
 {{\rm mor}}

\newcommand{\mstp}[3]% monomorphism version of \stp
 {\mbox{$#1\! : #2 \mono #3$}}

\newcommand{\Mu}%capital mu
 {{\rm M}}

\newcommand{\name}[1]% name of a relation
 {\mbox{$\ulcorner #1 \urcorner$}}

\newcommand{\names}[1]%\name used as subscript
 {\mbox{$\ulcorner$} #1 \mbox{$\urcorner$}}

\newcommand{\nml}% normal subgroup
 {\triangleleft}

\newcommand{\ob}% class of objects
 {{\rm ob}}

\newcommand{\op}% opposite category
 {^{\rm op}}

\newcommand{\pepi}% partial epimorphism
 {\rightharpoondown\kern-0.9em\rightharpoondown}

\newcommand{\pmap}% partial map arrow
 {\rightharpoondown}

\newcommand{\Pos}% positivization of a coherent category
 {{\bf Pos}}

\newcommand{\prarr}% parallel pair of arrows
 {\rightrightarrows}

\newcommand{\princfil}[1]% principal filter
 {\mbox{$\uparrow\!(#1)$}}

\newcommand{\princid}[1]% principal ideal
 {\mbox{$\downarrow\!(#1)$}}

\newcommand{\prstp}[3]% parallel-pair version of \stp
 {\mbox{$#1\! : #2 \prarr #3$}}

\newcommand{\pstp}[3]% partial-map version of \stp
 {\mbox{$#1\! : #2 \pmap #3$}}

\newcommand{\relarr}% relation-type arrow
 {\looparrowright}

\newcommand{\rlim}% right limit, i.e. colimit
 {{\mbox{$\lower.95ex\hbox{{\rm lim}}$}\atop{\scriptstyle
{\rightarrow}}}{}}

\newcommand{\rlimd}% \displaymath version of \rlim
 {\lower0.37ex\hbox{$\pile{\lim \\ {\scriptstyle
\rightarrow}}$}{}}

\newcommand{\rstp}[3]% relation version of \stp
 {\mbox{$#1\! : #2 \relarr #3$}}

\newcommand{\scn}% Sierpinski cone
 {{\bf scn}}

\newcommand{\scnS}% ditto relative to $\cal S$
 {{\bf scn}_{\cal S}}

\newcommand{\semid}% semidirect product
 {\rtimes}

\newcommand{\Sep}% category of separated objects
 {{\bf Sep}}

\newcommand{\sep}% category of separated objects
 {{\bf sep}}

\newcommand{\Set}% category of sets
 {{\bf Set }}

\newcommand{\Sh}% category of sheaves
 {{\bf Sh}}

\newcommand{\ShE}% sheaves relative to $\cal E$
 {{\bf Sh}_{\cal E}}

\newcommand{\ShS}% ditto for $\cal S$
 {{\bf Sh}_{\cal S}}

\newcommand{\sh}% category of sheaves
 {{\bf sh}}

\newcommand{\Simp}% the simplicial category
 {{\bf \Delta}}

\newcommand{\snd}% second projection
 {{\sf snd}}

\newcommand{\stg}[1]% string of #1
 {\vec{#1}}

\newcommand{\stp}[3]% source-target predicate
 {\mbox{$#1\! : #2 \to #3$}}

\newcommand{\Sub}% subobject lattice
 {{\rm Sub}}

\newcommand{\SUB}% indexed category of subobjects
 {{\mathbb S}{\sf ub}}

\newcommand{\tbel}% totally below
 {\prec\!\prec}

\newcommand{\tic}[2]%term-in-context, etc.
 {\mbox{$#1\!.\!#2$}}

\newcommand{\tp}% is of type
 {\!:}

\newcommand{\tps}% subscript version of \tp
 {:}

\newcommand{\tsub}% truncated subtraction
 {\pile{\lower0.5ex\hbox{.} \\ -}}

\newcommand{\wavy}% wavy arrow
 {\leadsto}

\newcommand{\wavydown}% wavy downarrow
 {\,{\mbox{\raise.2ex\hbox{\hbox{$\wr$}
\kern-.73em{\lower.5ex\hbox{$\scriptstyle{\vee}$}}}}}\,}

\newcommand{\wbel}% way-below relation
 {\lt\!\lt}

\newcommand{\wstp}[3]% wavy version of \stp
 {\mbox{$#1\!: #2 \wavy #3$}}

\newcommand{\fu}[2]
{[#1,#2]}

\newcommand{\st}[2]% source-target predicate
 {\mbox{$#1 \to #2$}}

\section{Introduction}

In a letter to Michael Makkai of 1989, Ieke Moerdijk proved the following result:

Let $\Sigma$ be a signature, and let $\Sigma$-\textbf{str}$({\cal E})$ denote the category of $\Sigma$-structures in a Grothendieck topos $\cal E$. Then a finitary first-order theory $\mathbb T$ over $\Sigma$ can be axiomatized by coherent sequents over $\Sigma$ if and only if 
 
(i) for any geometric morphism $f:{\cal F}\to {\cal E}$ between Grothendieck toposes, if $M\in \Sigma$-\textbf{str}$({\cal E})$ is a model of $\mathbb T$ then $f^{\ast}(M)$ is a model of $\mathbb T$;
 
(ii) for any surjective geometric morphism $f:{\cal F}\to {\cal E}$ between Grothendi-\\eck toposes and any $M\in \Sigma$-\textbf{str}$({\cal E})$, if $f^{\ast}(M)$ is a model of $\mathbb T$ then $M$ is a model of $\mathbb T$.

His proof of this result involved model-theoretic as well as topos-theoretic arguments, and heavily relied on the compactness theorem. %In fact, he first showed, by using arguments along the lines of the proof of the Tarski-Los theorem, that a finitary first-order theory is coherent if and only if its category of models is closed in $\Sigma$-\textbf{str}$({\Set })$ under filtered colimits, and then proved the equivalence of this latter condition with the conjunction of properties i) and ii) above.

In the same letter, Moerdijk asked for a proof of his conjecture that this result could be extended to the infinitary context i.e. that the version of it obtained by replacing `finitary first-order' with `infinitary first-order' and `coherent' by `geometric' also hold. This question remained unanswered for the past twenty years; in fact, the difficulty lies in the fact that, since in the infinitary context one can no longer rely on the compactness theorem, one cannot hope to prove the conjecture by extending the argument given in the finitary case.  

In this paper, we prove the conjecture by adopting the point of view of classifying toposes. We start by establishing some facts that will be useful for our analysis; then, in the third section, we prove our main theorem giving a semantic characterization of the classes of structures which arise as the collection of models in Grothendieck toposes of a geometric theory. In the last section, we derive Moerdijk's conjecture as an application of our criterion in the case of the class of models of an infinitary first-order theory, and we show that (a stronger version of) Moerdijk's result also follows as a consequence of our theorem. 

Before proceeding further, I would like to express my gratitude to Ieke Moerdijk for bringing my attention to his conjecture at a recent conference; it is also a pleasure to thank him, as well as Peter Johnstone, for their useful remarks on a preliminary version of this paper.

% prompting me to think about his conjecture

\section{Jointly surjective families of\\ geometric morphisms}

Recall from \cite{El} that a geometric morphism of (elementary) toposes is surjective if its inverse image functor is conservative i.e. it is faithful and reflects isomorphisms; more generally, a family $\{f_{i}:{\cal E}_{i}\to {\cal E} \textrm{ | } i\in I\}$ of geometric morphisms with common codomain is said to be jointly surjective if and only if the inverse image functors $f_{i}^{\ast}$ are jointly conservative. 

Note that if ${\cal C}$ and $\cal D$ are categories with equalizers and $F:{\cal C}\to {\cal D}$ is a functor preserving equalizers then $F$ is conservative if and only if it reflects isomorphisms; indeed, two arrows with common domain and codomain are equal if and only if their equalizer is an isomorphism. In particular, a family of geometric morphisms is jointly surjective if and only if the family formed by their inverse image functors jointly reflects isomorphisms.   
 
Given a collection $\{{\cal E}_{i} \hookrightarrow {\cal E} \textrm{ | } i\in I \}$ of subtoposes of a given elementary topos $\cal E$, we denote by $\mathbin{\mathop{\textrm{\huge $\cup$}}\limits_{i\in I}} {\cal E}_{i} \hookrightarrow {\cal E}$ the smallest subtopos of $\cal E$ containing all the ${\cal E}_{i}$, provided that it exists; recall from \cite{OC6} that if $\cal E$ is a Grothendieck topos then there is only a \emph{set} of (equivalence classes of) subtoposes of $\cal E$, and arbitrary unions of subtoposes always exist.

The following lemma gives a characterization of jointly surjective families of geometric morphisms.

\begin{lemma}\label{preslemma}
Let $\{f_{i}:{\cal E}_{i}\to {\cal E} \textrm{ | } i\in I\}$ be a family of geometric morphisms of elementary toposes with common codomain $\cal E$. Then $\{f_{i}:{\cal E}_{i}\to {\cal E} \textrm{ | } i\in I\}$ is jointly surjective if and only if ${\cal E}=\mathbin{\mathop{\textrm{\huge $\cup$}}\limits_{i\in I}} {\cal E}_{i}'$, where for each $i\in I$, ${\cal E}_{i}\epi {\cal E}_{i}' \hookrightarrow {\cal E}$ is the surjection-inclusion factorization of $f_{i}$.
\end{lemma}

\begin{proofs}

It is clear that $\{f_{i}:{\cal E}_{i}\to {\cal E} \textrm{ | } i\in I\}$  is jointly surjective if and only if the family $\{{\cal E}_{i}' \hookrightarrow {\cal E} \textrm{ | } i\in I\}$ of subtoposes of $\cal E$ is jointly surjective. For any $i\in I$, let $j_{i}$ denote the local operator on $\cal E$ corresponding to the subtopos ${\cal E}_{i}'$ of $\cal E$ and let $a_{j_{i}}:{\cal E}\to {\cal E}_{i}'$ be the corresponding associated sheaf functor. 

Let us suppose that $\{{\cal E}_{i}'\hookrightarrow {\cal E} \textrm{ | } i\in I\}$ is jointly surjective; we want to prove that ${\cal E}=\mathbin{\mathop{\textrm{\huge $\cup$}}\limits_{i\in I}} {\cal E}_{i}'$. Given a local operator $j$ on $\cal E$ which is smaller than each of the $j_{i}$, we want to prove that $j$ is the smallest local operator on $\cal E$. Now, for any arrow $f$ in $\cal E$, if $a_{j}(f)$ is an isomorphism then $a_{j_{i}}(f)$ is an isomorphism for each $i$, and hence, by our hypothesis, $f$ is an isomorphism; this proves our claim. 

Conversely, let us suppose that ${\cal E}=\mathbin{\mathop{\textrm{\huge $\cup$}}\limits_{i\in I}} {\cal E}_{i}'$; we have to prove that $\{{\cal E}_{i}'\hookrightarrow {\cal E} \textrm{ | } i\in I\}$ is jointly surjective i.e. for any arrow $f$ in $\cal E$, if $a_{j_{i}}(f)$ is an isomorphism  for every $i\in I$ then $f$ is an isomorphism. Now, for a fixed arrow $f$ in $\cal E$, consider the smallest local operator $k$ on $\cal E$ such that the corresponding associated sheaf functor $a_{k}$ sends $f$ to an isomorphism (cfr. Example A4.5.14(c) \cite{El}). By our hypothesis, $k\leq j_{i}$ for each $i$ and hence $k$ is the smallest local operator, which implies that $f$ is an isomorphism, as required.        
\end{proofs} 

\begin{rmk}\label{rem}
\emph{If all the toposes in the statement of the lemma are Grothendi-\\eck toposes and $I$ is a set then the lemma admits the following $2$-categorical interpretation. Recall that, for any set-indexed collection $\{{\cal E}_{i} \textrm{ | } i\in I\}$ of Grothendieck toposes, there exists the coproduct (Grothendieck) topos $\mathbin{\mathop{\textrm{$\coprod$}}\limits_{i\in I}} {\cal E}_{i}$. Now, it is immediate to see, by using the arguments in the proof of the lemma, that, given a family $\{f_{i}:{\cal E}_{i}\to {\cal E} \textrm{ | } i\in I\}$ of geometric morphisms with common codomain, the surjection-inclusion factorization of the induced coproduct map $f:\mathbin{\mathop{\textrm{$\coprod$}}\limits_{i\in I}} {\cal E}_{i}\to {\cal E}$ is given by its factorization through the inclusion $\mathbin{\mathop{\textrm{\huge $\cup$}}\limits_{i\in I}} {\cal E}_{i}'\hookrightarrow {\cal E}$; in particular, $\{f_{i}:{\cal E}_{i}\to {\cal E} \textrm{ | } i\in I\}$ is jointly surjective if and only if $f$ is surjective.}   
\end{rmk}

\section{The characterization theorem}

All the toposes in this section will be Grothendieck toposes.

Let $\Sigma$ be a signature. Let us denote by ${\mathbb O}_{\Sigma}$ the empty (geometric) theory over $\Sigma$ and by $\Set[{\mathbb O}_{\Sigma}]$ its classifying topos. Note that the ${\mathbb O}_{\Sigma}$-models in any Grothendieck topos $\cal E$ are precisely the $\Sigma$-structures in $\cal E$. Thus, for any Grothendieck topos $\cal E$, geometric morphisms ${\cal E} \to \Set[{\mathbb O}_{\Sigma}]$ correspond to $\Sigma$-structures in $\cal E$; the geometric morphism corresponding to a $\Sigma$-structure $M$ will be denoted by $f_{M}$ (note that if $U$ is a universal model of ${\mathbb O}_{\Sigma}$ in $\Set[{\mathbb O}_{\Sigma}]$ then $M\cong f_{M}^{\ast}(U)$).

Let us denote by $\Sigma$-\textbf{str}$({\cal E})$ the category of $\Sigma$-structures in a topos $\cal E$, as in the introduction above.

\begin{theorem}\label{teofond}

Let $\Sigma$ be a signature and $\cal S$ be a collection of $\Sigma$-structures in Grothendieck toposes closed under isomorphisms of structures. Then $\cal S$ is the collection of all models in Grothendieck toposes of a geometric theory over $\Sigma$ if and only if it satisfies the following two conditions:

(i) for any geometric morphism $f:{\cal F}\to {\cal E}$, if $M\in \Sigma$-\textbf{str}$({\cal E})$ is in $\cal S$ then $f^{\ast}(M)$ is in $\cal S$;

(ii) for any (set-indexed) jointly surjective family $\{f_{i}:{\cal E}_{i}\to {\cal E} \textrm{ | } i\in I\}$ of geometric morphisms and any $\Sigma$-structure $M$ in $\cal E$, if $f_{i}^{\ast}(M)$ is in $\cal S$ for every $i\in I$ then $M$ is in $\cal S$. 
\end{theorem}

\begin{proofs}
The `only if' part of the theorem is well-known. Let us prove the `if' part. Let us consider the collection of geometric morphisms to $\Set[{\mathbb O}_{\Sigma}]$ of the form $f_{M}$ for $M$ in $\cal S$; let ${\cal E}\hookrightarrow \Set[{\mathbb O}_{\Sigma}]$ be the subtopos of $\Set[{\mathbb O}_{\Sigma}]$ given by the union of all the subtoposes of $\Set[{\mathbb O}_{\Sigma}]$ arising as the inclusion parts of the surjection-inclusion factorizations of these geometric morphisms, and let $a:\Set[{\mathbb O}_{\Sigma}]\to {\cal E}$ be the corresponding associated sheaf functor. We know from \cite{OC6} (Theorem 3.6) that the subtopos ${\cal E}\hookrightarrow \Set[{\mathbb O}_{\Sigma}]$ of $\Set[{\mathbb O}_{\Sigma}]$ corresponds to a (unique up to syntactic equivalence) geometric quotient ${\mathbb T}$ of ${\mathbb O}_{\Sigma}$ such that if $U_{{\mathbb O}_{\Sigma}}$ is a universal model of ${\mathbb O}_{\Sigma}$ in $\Set[{\mathbb O}_{\Sigma}]$ then $U_{\mathbb T}:=a(U_{{\mathbb O}_{\Sigma}})$ is a universal model of $\mathbb T$ in $\cal E$. We will show that $\mathbb T$ axiomatizes our class of structures $\cal S$.

Let $M\in \Sigma$-\textbf{str}$({\cal E}_{M})$ be a structure in $\cal S$. The subtopos ${\cal E}_{M}'\hookrightarrow \Set[{\mathbb O}_{\Sigma}]$ arising in the surjection-inclusion factorization of $f_{M}:{\cal E}_{M}\to \Set[{\mathbb O}_{\Sigma}]$ factors as the inclusion ${\cal E}\hookrightarrow \Set[{\mathbb O}_{\Sigma}]$ composed with the canonical inclusion $l_{M}:{\cal E}_{M}'\hookrightarrow {\cal E}$. Now, if we compose this latter inclusion with the surjection part of the surjection-inclusion factorization of $f_{M}$, we obtain a geometric morphism $h_{M}:{\cal E}_{M}\to {\cal E}$ such that the composite of ${\cal E}\hookrightarrow \Set[{\mathbb O}_{\Sigma}]$ with $h_{M}$ is equal to $f_{M}$. But $M\cong f_{M}^{\ast}(U_{{\mathbb O}_{\Sigma}})$, from which it follows that $h_{M}^{\ast}(U_{\mathbb T})\cong M$, and hence that $M$ is a model of $\mathbb T$. This shows that every structure in $\cal S$ is a model of $\mathbb T$. To prove the converse, we note that, by Lemma \ref{preslemma}, the family of geometric morphisms $h_{M}$ for $M$ in $\cal S$ is jointly surjective; hence, under assumption (ii), $U_{\mathbb T}$ lies in $\cal S$. Now, since (by the universal property of the classifying topos $\cal E$ of $\mathbb T$) every model $N$ of $\mathbb T$ in a Grothendieck topos $\cal F$ is of the form $g^{\ast}(U_{{\mathbb T}})$ for some geometric morphism $g:{\cal F}\to {\cal E}$, condition (i) implies that any $\mathbb T$-model in a Grothendieck topos lies in $\cal S$. This concludes the proof of the theorem.

Note in passing that $\mathbb T$ can be described as the collection of all the geometric sequents over $\Sigma$ which are valid in every structure $M$ of $\cal S$ (cfr. also Theorem 9.1 \cite{OC6}). 
\end{proofs} 
 
It is natural to wonder if one can suppose the set $I$ in the statement of the theorem to be a singleton without loss of generality; in fact, we now show that this is not possible.

Given a class $\cal S$ of $\Sigma$-structures in Grothendieck toposes, we can explicitly describe the smallest class $\tilde{{\cal S}}$ of $\Sigma$-structures containing $\cal S$ which is closed under (i) and the version of (ii) obtained by requiring $I$ to have cardinality $1$. Indeed, with the notation used in the proof of the theorem, consider, for any $M$ in $\cal S$, the structure $\tilde{M}=i_{M}^{\ast}(U_{{\mathbb O}_{\Sigma}})$ where $i_{M}:{\cal E}_{M}'\to \Set[{\mathbb O}_{\Sigma}]$ is the inclusion part of the surjection-inclusion factorization of $f_{M}$; then $\tilde{{\cal S}}$ is equal to the collection $\cal R$ of all the $\Sigma$-structures of the form $g^{\ast}(\tilde{M})$ for some geometric morphism $g$. To prove this, we argue as follows. Clearly, $\cal R$ is contained in $\tilde{{\cal S}}$ and is closed under (i), so it remains to prove that it is closed under the version of (ii) obtained by requiring $I$ to have cardinality $1$. Let $N$ be a $\Sigma$-structure in a Grothendieck topos $\cal F$ and $p:{\cal G}\to {\cal F}$ be a surjective geometric morphism such that $p^{\ast}(N)$ is in $\cal R$; we want to prove that $N$ is in $\cal R$. Since $p^{\ast}(N)$ is in $\cal R$, there exists a $\Sigma$-structure $M$ in $\cal S$ such that $p^{\ast}(N)=g^{\ast}(\tilde{M})$ for some geometric morphism $g:{\cal G} \to {\cal E}_{M}'$. Then, by the universal property of the classifying topos for ${\mathbb O}_{\Sigma}$, the geometric morphisms $i_{M}\circ g$ and $f_{N}\circ p$ are isomorphic. Let ${\cal G}\stackrel{p_{g}}{\epi} {\cal U} \stackrel{g'}{\mono} {\cal E}_{M}'$ and ${\cal F} \stackrel{p_{f_{N}}}{\epi} {\cal F}' \stackrel{f_{N}'}{\mono} \Set[{\mathbb O}_{\Sigma}]$ be respectively the surjection-inclusion factorization of $g$ and of $f_{N}$; then ${\cal G}\stackrel{p_{g}}{\epi} {\cal U} \stackrel{i_{M}\circ g'}{\mono} \Set[{\mathbb O}_{\Sigma}]$ and ${\cal G} \stackrel{p_{f_{N}}\circ p}{\epi} {\cal F}' \stackrel{f_{N}'}{\mono} \Set[{\mathbb O}_{\Sigma}]$ are respectively the surjection-inclusion factorization of $i_{M}\circ g$ and of $f_{N}\circ p$. Then, by the uniqueness (up to equivalence) of the surjection-inclusion factorization of a geometric morphism, the geometric morphisms $i_{M}\circ g'$ and $f_{N}'$ are isomorphic, from which it follows that $N$ is in $\cal R$. This completes the proof of the equality ${\cal R}=\tilde{{\cal S}}$. 

Let us now show that it is not true in general that $\tilde{{\cal S}}$ is axiomatized by a geometric theory over $\Sigma$. For a counterexample, take two subtoposes $i_{1}:{\cal E}_{1}\hookrightarrow \Set[{\mathbb O}_{\Sigma}]$ and $i_{2}:{\cal E}_{2}\hookrightarrow \Set[{\mathbb O}_{\Sigma}]$ of $\Set[{\mathbb O}_{\Sigma}]$ which are not contained in each other, and take $\cal S$ to consist of the two models $M_{1}:=i_{1}^{\ast}(U_{{\mathbb O}_{\Sigma}})$ and $M_{2}:=i_{2}^{\ast}(U_{{\mathbb O}_{\Sigma}})$; if $\tilde{\cal S}$ were axiomatized by a geometric theory over $\Sigma$ then, by Lemma \ref{preslemma}, the $\Sigma$-structure $i^{\ast}(U_{{\mathbb O}_{\Sigma}})$, where $i:{\cal E}_{1}\cup {\cal E}_{2}\hookrightarrow \Set[{\mathbb O}_{\Sigma}]$ is the union of the subtoposes ${\cal E}_{1}$ and ${\cal E}_{2}$ of $\Set[{\mathbb O}_{\Sigma}]$, would lie in $\tilde{{\cal S}}$, and we can show this to be impossible. Indeed, if $i^{\ast}(U_{{\mathbb O}_{\Sigma}})$ were in $\tilde{{\cal S}}$ then there would be a geometric morphism $g:{\cal E}_{1}\cup {\cal E}_{2}\to {\cal E}_{1}$ (or $g:{\cal E}_{1}\cup {\cal E}_{2}\to {\cal E}_{2}$) such that $i\cong i_{1} \circ g$ (or $i\cong i_{2}\circ g$); but the existence and uniqueness of the surjection-inclusion factorizations of a geometric morphism ensure that $g$ is an equivalence, which contradicts our assumption that ${\cal E}_{1}$ and ${\cal E}_{2}$ be not contained in each other.
      
\begin{rmk}\label{rmk2}
\emph{In view of Remark \ref{rem}, condition (ii) in the statement of the theorem can be rephrased as follows:}

\emph{(i) for any surjective geometric morphism $f:{\cal F}\to {\cal E}$ and any $M\in \Sigma$-\textbf{str}$({\cal E})$, if $f^{\ast}(M)$ is in $\cal S$ then $M$ is in $\cal S$;}

\emph{(ii) for any set-indexed family $\{M_{i} \textrm{ | } i \in I\}$ of structures in toposes ${\cal E}_{i}$ all of which are in $\cal S$, the structure in the coproduct topos $\mathbin{\mathop{\textrm{$\coprod$}}\limits_{i\in I}} {\cal E}_{i}$ whose $i$th coordinate is $M_{i}$ is also in $\cal S$.}
\end{rmk}

\section{Applications}

Let $\mathbb T$ be an infinitary first-order theory over a given signature $\Sigma$ and ${\cal S}_{\mathbb T}$ be the collection of its models inside Grothendieck toposes. Clearly, ${\cal S}_{\mathbb T}$ satisfies condition (ii) of Remark \ref{rmk2}, so it is axiomatizable by geometric sequents over $\Sigma$ if and only if it satisfies condition (i) of Theorem \ref{teofond} and condition (i) of Remark \ref{rmk2}. Note that, by Proposition D1.3.2 \cite{El2} and Corollary 3.4 \cite{BJ}, two infinitary first-order theories over the same signature are deductively equivalent (relative to the intuitionistic proof system of infinitary first-order logic of section D1.3 \cite{El2}) if and only if they have the same models in every Grothendieck topos. Hence we have proved the conjecture by Moerdijk mentioned in the introduction of this paper.

If $\mathbb T$ is a finitary first-order theory satisfying the conditions of the characterization theorem, it is natural to wonder whether the geometric theory axiomatizing $\mathbb T$ provided by the theorem is in fact coherent. As recalled in the introduction above, it was already proved by Moerdijk in his letter that if $\mathbb T$ satisfies condition (i) of Theorem \ref{teofond} and condition (i) of Remark \ref{rmk2} then $\mathbb T$ is axiomatizable over $\Set$ by coherent sequents over its signature. In fact, it will follow directly from our theorem that this is true not only over $\Set$ but over every Grothendieck topos, once we have shown that if $\mathbb T$ is a finitary first-order theory over a signature $\Sigma$ and ${\mathbb T}'$ is a geometric theory over $\Sigma$ having the same models in Grothendieck toposes as ${\mathbb T}$ then ${\mathbb T}'$ is coherent. To prove this, we argue as follows.     

By using Theorem 3.5 \cite{OC8}, we are reduced to verify that for any coherent formula $\{\vec{x}. \phi\}$ over $\Sigma$, for any family $\{\psi_{i}(\vec{x}) \textrm{ | } i\in I\}$ of coherent (equivalently, geometric) formulae in the same context, if $\phi \vdash_{\vec{x}} \mathbin{\mathop{\textrm{\huge $\vee$}}\limits_{i\in I}} \psi_{i}$ is provable in ${\mathbb T}'$ (using geometric logic) then $\phi \vdash_{\vec{x}} \mathbin{\mathop{\textrm{\huge $\vee$}}\limits_{i\in I'}} \psi_{i}$ is provable in ${\mathbb T}'$ (using geometric logic) for some finite subset $I'$ of $I$.

We can suppose, without loss of generality, $\vec{x}$ to be the empty string; indeed, if $\vec{c}$ is a string of new constants of the same length and type as $\vec{x}$, a geometric sequent $\chi \vdash_{\vec{x}} \xi$ over $\Sigma$ is provable in ${\mathbb T}'$ if and only if the sequent $\chi[\vec{c}\slash \vec{x}] \vdash_{[]} \xi[\vec{c}\slash \vec{x}]$ is provable in ${\mathbb T}'$, regarded as a theory over the signature $\Sigma \cup \{\vec{c}\}$.

If $\phi \vdash_{\vec{x}} \mathbin{\mathop{\textrm{\huge $\vee$}}\limits_{i\in I}} \psi_{i}$ is provable in ${\mathbb T}'$ then every model in $\Set$ of the theory ${\mathbb T} \cup \{\neg \psi_{i} \textrm{ | } i\in I\}$ is a model of $\neg \phi$. Since the theory ${\mathbb T} \cup \{\neg \psi_{i} \textrm{ | } i\in I\}$ is finitary first-order, this condition is equivalent to saying that $\neg \phi$ is provable in the theory ${\mathbb T} \cup \{\neg \psi_{i} \textrm{ | } i\in I\}$ (using classical finitary first-order logic), from which it follows, by the finiteness theorem in classical Model Theory, that $\neg \phi$ is provable in ${\mathbb T} \cup \{\neg \psi_{i}  \textrm{ | } i\in I'\}$ for some finite subset $I'$ of $I$ i.e. $\phi \vdash_{[]} \mathbin{\mathop{\textrm{\huge $\vee$}}\limits_{i\in I'}} \psi_{i}$ is provable in $\mathbb T$ (using classical finitary first-order logic). Thus $\phi \vdash_{\vec{x}} \mathbin{\mathop{\textrm{\huge $\vee$}}\limits_{i\in I'}} \psi_{i}$ is valid in every model of $\mathbb T$ (equivalently, of ${\mathbb T}'$) in Boolean Grothendieck toposes and hence, by Proposition D3.1.16 \cite{El2}, it is provable in ${\mathbb T}'$, as required.

\end{document}